\documentclass[11pt]{article}

\usepackage[dvips]{graphicx}
\usepackage{amssymb}
\usepackage{amsmath}
\usepackage{amsthm}
\usepackage{latexsym}
\usepackage[english]{babel}
\usepackage{appendix}

\usepackage[usenames,dvipsnames]{color}

\newtheorem{theorem}{Theorem}[section]
\newtheorem{definition}{Definition}[section]
\newtheorem{lemma}{Lemma}[section]
\newtheorem{proposition}{Proposition}[section]

\newtheorem{corollary}{Corollary}[section]

\def \C {\mathbb{C}}
\def \R {\mathbb{R}}
\def \H {\mathbb{H}}
\def \N {\mathbb{N}}

\def \L {\mathcal{L}}

\def \P {\mathcal{P}}
\def \Cay {\mathcal{C}}

\def\Xint#1{\mathchoice
{\XXint\displaystyle\textstyle{#1}}%
{\XXint\textstyle\scriptstyle{#1}}%
{\XXint\scriptstyle\scriptscriptstyle{#1}}%
{\XXint\scriptscriptstyle\scriptscriptstyle{#1}}%
\!\int}
\def\XXint#1#2#3{{\setbox0=\hbox{$#1{#2#3}{\int}$ }
\vcenter{\hbox{$#2#3$ }}\kern-.6\wd0}}

\def\dashint{\Xint-}

\newcommand{\ol}[1]{\overline{#1}}

\begin{document}

\title{Logarithmic Hardy-Littlewood-Sobolev Inequality on Pseudo-Einstein $3$-manifolds and the Logarithmic Robin Mass}

\author{Ali Maalaoui$^{(1)}$}
\addtocounter{footnote}{1}
\footnotetext{Department of mathematics and natural sciences, American University of Ras Al Khaimah, PO Box 10021, Ras Al Khaimah, UAE. E-mail address:
{\tt{ali.maalaoui@aurak.ac.ae}}}

\date{}
\maketitle

\vspace{5mm}

{\noindent\bf Abstract} {Given a three dimensional pseudo-Einstein CR manifold $(M,T^{1,0}M,\theta)$, we study the existence of a contact structure conformal to $\theta$ for which the logarithmic Hardy-Littlewood-Sobolev (LHLS) inequality holds. Our approach closely follows \cite{Ok1} in the Riemannian setting. For this purpose, we introduce the notion of Robin mass as the constant term appearing in the expansion of the Green's function of the $P'$-operator. We show that the LHLS inequality appears when we study the variation of the total mass under conformal change. Then we exhibit an Aubin type result guaranteeing the existence of a minimizer for the total mass which yields the classical LHLS inequality.}

\vspace{5mm}


\vspace{4mm}


\vspace{4mm}


\section{Introduction  and statement of the results}

The logarithmic Hardy-Littlewood-Sobolev inequality (LHLS) is one of the most fundamental inequalities in analysis since it appears as the borderline case of the classical Hardy-Littlewood Sobolev inequalities which in their turn are the dual of the classical Sobolev embeddings. We refer the reader for instance to \cite{CL,L} and the references therein. We recall that in the standard sphere $(S^{n},g_{0})$, this inequality reads as:
\begin{equation}\label{lhsc}
\frac{2}{n!}\int_{S^{n}}F\ln(F)\ dv_{g_{0}}-\int_{S^{n}}FA_{n}^{-1}F\ dv_{g_{0}}\geq 0,
\end{equation}
for all $F:S^{n}\to \R_{+}$ such that $\int_{S^{n}}F\ dv=1$ with $\int_{S^{n}}F\ln(F)\ dv<\infty$. Here $A_{n}$ is the Paneitz operator defined by its action on the spherical harmonics $Y_{k}$ by
$$A_{n}Y_{k}=k(k+1)\cdots (k+n-1)Y_{k}.$$
 The dual of $(\ref{lhsc})$ is the classical Beckner-Onofry inequality \cite{Beck,CL,CY}, that states that for $u\in H^{\frac{n}{2}}(S^{n})$, 
$$
\frac{1}{2n!}\dashint_{S^{n}}uA_{n}u\ dv_{g_{0}}+\dashint_{S^{n}}u\ dv_{g_{0}} -\ln\Big(\dashint_{S^{n}} e^{u}\ dv_{g_{0}}\Big)\geq 0.
$$

From a spectral point of view, the LHLS inequality appears in estimating the regularized spectral zeta function of the operator $A_{n}$ as proved in \cite{Mor1}: if $\int_{S^{n}}Fdv_{g_{0}}=1$, then
\begin{equation}\label{zeta}
\tilde{Z}(\tilde{g})-\tilde{Z}(g_{0})=\frac{2}{n!}\int_{S^{n}}F\ln(F)\ dv_{g_{0}}-\int_{S^{n}}FA_{n}^{-1}F\ dv_{g_{0}},
\end{equation}
where $\tilde{g}=F^{\frac{2}{n}}g_{0}$ and $\tilde{Z}$ is the regularized Zeta function of the operator $A_{n}$ which also can be replaced by $trace (A_{n}^{-1})$ as in \cite{S1,S2,Ok1}. This spectral property was then investigated in \cite{Ok1}, in the case of general Riemannian manifolds:
\begin{theorem}[\cite{Ok1}]
Let $\Gamma_{V}$ be a conformal a class of metrics on $M^{n}$ with a fixed volume $V$, then
$$ \inf_{g\in \Gamma_{V}}trace A_{n}^{-1}(M,V)\leq trace A_{n}^{-1}(S^{n},V),$$
where $A_{n}(M,V)$ is the critical GJMS operator on the manifold $M$ with volume $V$ (\cite{GJMS}). Moreover, if the inequality is strict, then the infimum is attained by a metric in $\Gamma_{V}$.
\end{theorem}
This result was proved by introducing a notion of mass for the Green's function of the critical GJMS operator (see \cite{Ok1,S1,S2}). Indeed, as in the case of the mass for the Yamabe-type problems, the Robin mass is the constant term appearing after the logarithmic singularity in the expansion of the Green's function.\\

In this work we will focus on the three dimensional CR setting. With this setting, there are fundamental differences compared to the Riemannian setting. In fact, one does not have a general Moser-Trudinger inequality unless the study is restricted to pluriharmonic function $\P$. In fact, the right substitute for the critical GJMS operator in this case, is the $P^{\prime}$-Paneitz type operator. For instance, in $S^{2n+1}$ this operator acts as
$$ P^{\prime}\sum_{j}( Y_{0,j}+Y_{j,0})=\sum_{j}\lambda_{j}(Y_{0,j}+Y_{j,0}),$$
where $\lambda_{j}=j(j+1)\cdots (j+n)$ and $Y_{0,j},Y_{j,0}$ form an $L^{2}$-orthonormal basis of pluriharmonic functions on $S^{2n+1}$. Moreover, as shown in \cite{Bran}, one has a the following Moser-Trudinger inequality 
$$\frac{1}{2(n+1)!}\dashint_{S^{2n+1}} FP^{\prime}F\ dv+\dashint_{S^{2n+1}}F\ dv-\ln \dashint_{S^{2n+1}}e^{F}\ dv\geq 0,$$
for $F\in \P\cap W^{2,2}(M)$. Its dual, can be stated as follows: for any $G:S^{2n+1}\to \R$ with $G\geq 0$, $G\in LlogL$ and $\dashint_{S^{2n+1}}G\ dv=1$, we have
$$\frac{(n+1)!}{2}\dashint_{S^{2n+1}}(G-1)P^{\prime-1}\tau(G-1)\ dv\leq \dashint_{S^{2n+1}}G\ln(G)\ dv.$$
In fact if one looks closely to the work \cite{Bran}, one sees that the operator in case is actually $\overline{P}^{\prime}:=\tau P^{\prime}$ where $\tau$ is the $L^{2}$-projection on $\P$.
In this paper we propose to study the notion of Robin mass in the three dimensional CR setting and relate it to the LHLS inequality. Indeed, given an embeddable pseudo-Einstein manifold $(M,T^{1,0},\theta)$, then the $\overline{P}^{\prime}$ operator is well defined, and its Green's function $G_{\theta}$ takes the form 
$$G_{\theta}(x,y)=-\gamma_{3}\ln(d_{\theta}(x,y))+O(1),$$
where $\gamma_{3}=\frac{1}{4\pi^{2}}$.

\begin{definition}
Given an embeddable pseudo-Einstein manifold $(M,T^{1,0}M,\theta)$, the CR-Robin mass is defined by
\begin{equation}
m_{\theta}(x)=\lim_{y\to x} G_{\theta}(x,y)+\gamma_{3}\ln(d_{\theta}(x,y)).
\end{equation}
where $d_{\theta}$ is the horizontal distance induced by the Levi form $L_{\theta}$.\\
The total mass of $(M,T^{1,0}M,\theta)$ is then defined by
 $$\mathcal{M}_{\theta}(M):=\int_{M}m_{\theta}\ dv_{\theta}.$$
\end{definition}
For example, an easy computation in the case of the standard sphere $(S^{3},T^{1,0}S^{3},\theta_{0})$ yields $$m_{\theta_{0}}=\frac{1}{8\pi^{2}}\ln(2).$$
We define the space $\L(M)$ by
$$\L(M):=L^{1}Log(L^{1})^{+}(M)=\{F:M\to \R^{+}; \int_{M}F\ln(F)\ \ dv_{\theta}<\infty\}.$$
On the Heisenberg group $\H$, we let 
$$\L_{c}(\H):=\{f:\H\to \R^{+}; \text{ f is measurable and compactly supported and } \int_{\H}f\ln(f)dx<\infty\}.$$
We fix now a pseudo-Einstein structure $(M,T^{1,0}M,\theta)$, such that $P'_{\theta}$ is non-negative and $\ker P'_{\theta}=\R$ and without loss of generality we can assume that
$$V:=\int_{M}\theta\wedge d\theta=\int_{S^{3}} \theta_{0}\wedge d\theta_{0}.$$
 We also let $\theta_{F}=F^{\frac{1}{2}}\theta$. Notice that $\theta_{F}$ induces a pseudo-Einstein structure if and only if $\ln(F)\in \mathcal{P}$. We set $\tau_{F}$ to be the orthogonal projection on $\P$ with respect to the $L^{2}$-inner product induced by $\theta_{F}$. We can then define the operators
$$A_{\theta}:=\tau P^{\prime}_{\theta}\tau \qquad\text{ and }\qquad A_{\theta_{F}}=\tau_{F}\Big( F^{-1}A_{\theta}\Big).$$
We will also write $V_{F}$ for the volume of $M$ with respect to $dv_{\theta_{F}}$, that is,
$$V_{F}=\int_{M}\theta_{F}\wedge d \theta_{F}=\int_{M}F\theta\wedge d\theta.$$
\textbf{Conventions:}\\
Under the assumption that $P'_{\theta}$ is non-negative and $\ker P'_{\theta}=\R$, we will make here a very important convention that will be carried along all the paper:\\
{\it Consider the operator $A_{\theta}$, then we will let the operator $A^{-1}_{\theta}$ act on all functions in $\mathcal{P}$ with the convention that 
$$A_{\theta}\circ A^{-1}_{\theta}\tau u=\tau u-\frac{1}{V}\int_{M}u \ dv_{\theta} \text{ and } A^{-1}_{\theta}1=0.$$
}
We will state below our main results and when there is no confusion, we will drop the dependence of the total mass on the manifold $M$. We have then the following result which is the CR version of $(\ref{zeta})$ proved in \cite{Mor1}.
\begin{theorem}\label{thrm1}
Consider an embeddable pseudo-Einstein manifold $(M,T^{1,0}M,\theta)$ such that $A_{\theta}$ is non-negative and $\ker A_{\theta}=\R$.\\
 If $m_{\theta}$ is constant and under the constraint of $V_{F}=V$, one has:
$$\mathcal{M}_{\theta_{F}}-\mathcal{M}_{\theta}=\frac{\gamma_{3}}{4}\int_{M}F\ln(F)\ dv_{\theta}-\frac{1}{V}\int_{M}FA_{\theta}^{-1}\tau F\ dv_{\theta},$$
for all $F\in \L(M)$. In particular, on the standard sphere $(S^{3},T^{1,0}S^{3},\theta_{0})$ one has
$$\mathcal{M}_{\theta_{F}}(S^{3})-\mathcal{M}_{\theta_{0}}(S^{3})=\frac{\gamma_{3}}{4}\int_{S^{3}}F\ln(F)\ dv_{\theta}-\frac{1}{V}\int_{S^{3}}FA_{\theta}^{-1}\tau F\ dv_{\theta}\geq 0.$$
with equality if and only if $F=|J_{k}|$ with $k \in Aut(S^{3})$, normalized to have volume $V$.
\end{theorem}
As we will show in Appendix A, if $\tilde{Z}_{\theta}(1)$ is the regularized Zeta function of $A_{\theta}$, then there exists a constant $c$ such that
$$\mathcal{M}_{\theta}=\tilde{Z}_{\theta}(1)+c.$$
Therefore, the previous Theorem can be reformulated as follow:
\begin{corollary}
Assume that $(M,T^{1,0}M,\theta)$ is an embeddable pseudo-Einstein manifold such that $A_{\theta}$ is non-negative and $\ker A_{\theta}=\R$. If $m_{\theta}$ is constant, then under the constraint $V_{F}=V$, we have
$$\tilde{Z}_{\theta_{F}}(1)-\tilde{Z}_{\theta}(1)=\frac{\gamma_{3}}{4}\int_{M}F\ln(F)\ dv_{\theta}-\frac{1}{V}\int_{M}FA_{\theta}^{-1}\tau F\ dv_{\theta},$$
for all $F\in \L(M)$. In particular, on the standard sphere $(S^{3},T^{1,0}S^{3},\theta_{0})$ one has
\begin{equation}\label{zet}
\tilde{Z}_{\theta_{F},S^{3}}(1)-\tilde{Z}_{\theta_{0},S^{3}}(1)=\frac{\gamma_{3}}{4}\int_{S^{3}}F\ln(F)\ dv_{\theta}-\frac{1}{V}\int_{S^{3}}FA_{\theta}^{-1}\tau F\ dv_{\theta}\geq 0.
\end{equation}
with equality if and only if $F=|J_{k}|$ with $k \in Aut(S^{3})$, normalized to have volume $V$.
\end{corollary}
A truncated version of $(\ref{zet})$ was proved in \cite[Proposition 3.5]{Bran}. Indeed, if $\lambda_{k}(\theta)$ denotes the $k^{th}$ eigenvalue of $A_{\theta}$, then the authors show that
$$\sum_{k=1}^{4}\frac{1}{\lambda_{k}(\theta_{F})}\geq\sum_{k=1}^{4}\frac{1}{\lambda_{k}(\theta_{0})}.$$

Next, we will deduce a result that can be seen as an Aubin type result as in \cite{Aub} for the Yamabe problem and \cite{JL} for the CR-Yamabe problem. 
\begin{theorem}\label{thrm2}
We define
$$\mathcal{M}([\theta],M):=\inf_{F\in \L(M);V_{F}=V}\mathcal{M}_{\theta_{F}}(M).$$
Then under the assumptions of Theorem \ref{thrm1}, we have
\begin{itemize}
\item[i)] $\mathcal{M}([\theta],M)\leq \mathcal{M}([\theta_{0}],S^{3}).$
\item [ii)] If $\mathcal{M}([\theta],M)< \mathcal{M}([\theta_{0}],S^{3})$, then the infimum is achieved.
\end{itemize}

Moreover, if $\mathcal{M}([\theta],M)$ is achieved by a function $F_{0}$, then the contact form $\theta_{F_{0}}$ has constant mass and the LHLS holds, i.e, for all $F\in \L(M)$, such that $V=V_{F}$, we have
$$\frac{\gamma_{3}}{4}\int_{M}F\ln(F)\ dv_{\theta_{F_{0}}}-\frac{1}{V}\int_{M}FA^{-1}_{\theta_{F_{0}}}\tau F\ dv_{\theta_{F_{0}}}\geq 0.$$
If in addition $m_{\theta} \in \mathcal{P}$, then $\theta_{F_{0}}$ is pseudo-Einstein.
\end{theorem}
Based on the work in \cite{Chan} and \cite{CaYa}, the assumptions that $M$ is embeddable and $P'_{\theta}$ is non-negative and $\ker P'_{\theta}=\R$, can be replaced by the non-negativity of the Paneitz operator $P_{\theta}$ and that the conformal class $[\theta]$ carries a pseudo-Einstein structure with non-negative Webster curvature but non-identically zero.\\
One is also hoping to have a positive mass type theorem as in \cite{CMY}, stating that if $\mathcal{M}([\theta],M)=\mathcal{M}([\theta_{0}],S^{3})$, then $(M,\theta)$ is CR-equivalent to the standard sphere $(S^{3},\theta_{0})$, but for now, this type of result is beyond the work done in this paper and it needs a more refined blow-up analysis of the functional $J(\cdot,M)$ defined below. \\

\textbf{Acknowledgement:} The author wants to express his thanks to Carlo Morpurgo for his helpful comments on an earlier version of this paper.

\section{Preliminaries and Setting}
In this section we survey the main quantities and properties that we will be using during our investigation.
\subsection{Pseudo-Hermitian geometry}
We will closely follow the notations in \cite{CaYa}. Let $M$ be a smooth, oriented three-dimensional manifold.  A CR structure on $M$ is a one-dimensional complex subbundle $T^{1,0}\subset T_{\C}M:= TM\otimes\C$ such that $T^{1,0}\cap T^{0,1}=\{0\}$ for $T^{0,1}:=\overline{T^{1,0}}$.  Let $H=Re T^{1,0}$ and let $J\colon H\to H$ be the almost complex structure defined by $J(Z+\bar Z)=i(Z-\bar Z)$, for all $Z\in T^{1,0}$. The condition that $T^{1,0}\cap T^{0,1}=\{0\}$ is equivalent to the existence of a contact form $\theta$ such that $\ker \theta =H$. We recall that a 1-form $\theta$ is said to be a contact form if $\theta\wedge d\theta$ is a volume form on $M$. Since $M$ is oriented, a contact form always exists, and is determined up to multiplication by a positive real-valued smooth function. We say that $(M,T^{1,0}M)$ is strictly pseudo-convex if the Levi form $d\theta(\cdot,J\cdot)$ on $H\otimes H$ is positive definite for some, and hence any, choice of contact form $\theta$. We shall always assume that our CR manifolds are strictly pseudo-convex.

Notice that in a CR-manifold, there is no canonical choice of the contact form $\theta$. A pseudohermitian manifold is a triple $(M,T^{1,0}M,\theta)$ consisting of a CR manifold and a contact form. The Reeb vector field $T$ is the vector field such that $\theta(T)=1$ and $d\theta(T,\cdot)=0$. The choice of $\theta$ induces a natural $L^{2}$-dot product $\langle\cdot,\cdot\rangle$, defined by
$$\langle f,g\rangle =\int_{M}f(x)g(x)\theta\wedge d\theta.$$

 A $(1,0)$-form is a section of $T_{\C}^\ast M$ which annihilates $T^{0,1}$.  An admissible coframe is a non-vanishing $(1,0)$-form $\theta^1$ in an open set $U\subset M$ such that $\theta^1(T)=0$.  Let $\theta^{\bar 1}:=\overline{\theta^1}$ be its conjugate.  Then $d\theta=ih_{1\bar 1}\theta^1\wedge\theta^{\bar 1}$ for some positive function $h_{1\bar 1}$.  The function $h_{1\bar 1}$ is equivalent to the Levi form.  We set $\{Z_1,Z_{\bar 1},T\}$ to the dual of $(\theta^{1},\theta^{\bar 1},\theta)$. The geometric structure of a CR manifold is determined by the  connection form $\omega_1{}^1$ and the torsion form $\tau_1=A_{11}\theta^1$ defined in an admissible coframe $\theta^1$ and is uniquely determined by
\begin{equation*}
\left\{\begin{array}{ll}
d\theta^1 = \theta^1\wedge\omega_1{}^1 + \theta\wedge\tau^1, \\
\omega_{1\bar 1} + \omega_{\bar 11} = dh_{1\bar 1},
\end{array}
\right.
\end{equation*}
where we use $h_{1\bar 1}$ to raise and lower indices.  The connection forms determine the pseudohermitian connection $\nabla$, also called the Tanaka-Webster connection, by
\[ \nabla Z_1 := \omega_1{}^1\otimes Z_1. \] 
The scalar curvature $R$ of $\theta$, also called the Webster curvature, is given by the expression
 \[ d\omega_1{}^1 = R\theta^1\wedge\theta^{\bar 1} \mod\theta . \]
\begin{definition}
A real-valued function $w\in C^\infty(M)$ is CR pluriharmonic if locally $w=Re f$ for some complex-valued function $f\in C^\infty(M,\C)$ satisfying $Z_{\bar 1}f=0$. 
\end{definition}
Equivalently, \cite{Lee}, $w$ is a CR pluriharmonic function if
\[ P_{3} w:=\nabla_1\nabla_1\nabla^1 w + iA_{11}\nabla^1 w = 0, \]
for $\nabla_1:=\nabla_{Z_1}$. We denote by $\mathcal{P}$ the space of all CR pluriharmonic functions and $\tau: L^{2}(M)\to L^{2}(M)\cap \mathcal{P}$ be the orthogonal projection on the space of pluriharmonic functions. If $S:L^{2}(M)\to \ker \bar\partial_{b}$ denotes the Szeg\"{o} kernel, then 
$$\tau=S+\bar S+\mathcal{F},$$
where $\mathcal{F}F$ is a smoothing kernel as shown in \cite{H}. In particular, one has that $\tau$ is a bounded operator from $W^{k,p}(M)\to W^{k,p}(M)$ for $1<p<\infty$ and $k\in \N$ (see \cite{PS}). In fact, this last property can be directly deduced from the work \cite{H}, since in the author provides an expansion of the kernel of $\tau$ that we will still denote it by $\tau$:
\begin{theorem}[\cite{H}]\label{texp}
Assume that $(M,T^{1.0}M)$ is a compact embeddable strongly pseudo-convex CR manifold, then there exist $F_{1},G_{1}\in C^{\infty}(M\times M)$ such that
$$\tau(x,y)=2\text{Re}\Big(F_{1}(-i\varphi(x,y))^{-2}+G_{1}\ln(-i\varphi(x,y))\Big),$$
with $$F_{1}=a_{0}(x,y)+a_{1}(x,y)(-i\varphi(x,y))+f_{1}(x,y)(-i\varphi(x,y))^{2},$$ 
where $f_{1}\in C^{\infty}(M\times M)$ and $\varphi(x,y)$ has the following expansion in local coordinates near $x_{0}\in M$ and $x=x_{3}+iz$ and $y=y_{3}+iw$,
\begin{align}
\varphi(x,y)=&-x_{3}+y_{3}+i|z-w|^{2}+\Big(i(\bar{z}w-z\bar{w})+c(-zx_{3}+wy_{3})+\bar{c}(-\bar{z}x_{3}+\bar{w}y_{3})\Big)\notag\\
&+|x_{3}-y_{3}|f(x,y)+O(|(x,y)|^{3}),\notag
\end{align}
and $f$ is a smooth real function such that $f(0,0)=0$.
\end{theorem}
In particular, one can check that the first term of the expansion of $\tau$ coincides with the real part of the Szeg\"{o} projection in $\H$.

The Paneitz operator $P_\theta$ is the differential operator
\begin{align*}
P_\theta(w) & := 4\text{div}(P_{3}w) \\
& = \Delta_b^2w + T^2 - 4\text{Im} \nabla^1\left(A_{11}\nabla^1f\right),
\end{align*}
for $\Delta_b:=\nabla^1\nabla_1+\nabla^{\bar 1}\nabla_{\bar 1}$ the subLaplacian.  In particular, $\mathcal{P}\subset\ker P_\theta$. Hence, $\ker P_{\theta}$ is infinite dimensional. For a thorough study of the analytical properties of $P_{\theta}$ and its kernel, we refer the reader to \cite{H,CCY,CaYa1}. The main property of the Paneitz operator $P_{\theta}$ is that it is CR covariant \cite{Hir2}. That is, if $\hat\theta=e^w\theta$, then $e^{2w} P_{\hat{\theta}}=P_\theta$.
\begin{definition}
\label{pprime}
Let $(M^3,T^{1,0}M,\theta)$ be a pseudohermitian manifold.  The Paneitz type operator $P_{\theta}^\prime\colon\mathcal{P}\to C^\infty(M)$ is defined by
\begin{align}
P_\theta^\prime f & = 4\Delta_b^2 f - 8 \textnormal{Im}\left(\nabla^\alpha(A_{\alpha\beta}\nabla^\beta f)\right) - 4 \textnormal{Re}\left(\nabla^\alpha(R\nabla_\alpha f)\right) \notag\\
& \quad + \frac{8}{3}\textnormal{Re} (\nabla_\alpha R - i\nabla^\beta A_{\alpha\beta})\nabla^\alpha f - \frac{4}{3}f\nabla^\alpha( \nabla_\alpha R - i\nabla^\beta A_{\alpha\beta})\label{pprime1}
\end{align}
for $f\in\mathcal{P}$.
\end{definition}
The main property of the operator $P_{\theta}^\prime$ is its "almost" conformal covariance as shown in \cite{BG,CaYa}. That is if  $(M,T^{1,0}M,\theta)$ is a pseudohermitian manifold, $w\in C^\infty(M)$, and we set $\hat\theta=e^w\theta$, then
\begin{equation}
 \label{pprime}
 e^{2w} P_{\hat\theta}^\prime(u) = P_\theta^\prime(u) + P_\theta\left(uw\right)
\end{equation}
for all $u\in\mathcal{P}$.  In particular, since $P_\theta$ is self-adjoint and $\mathcal{P}\subset \ker P_{\theta}$, we have that the operator $P^\prime$ is conformally covariant, mod $\mathcal{P}^\perp$.
\begin{definition}
A pseudohermitian manifold $(M,T^{1,0}M,\theta)$ is pseudo-Einstein if $$\nabla_\alpha R - i\nabla^\beta A_{\alpha\beta}=0.$$
\end{definition}
 Moreover, if $\theta$ induces a pseudo-Einstein structure then $e^{u}\theta$ is pseudo-Einstein if and only if $u\in \mathcal{P}$. The definition above was stated in \cite{CaYa}, but it was implicitly mentionned in \cite{Hir2}. In particular, if $(M^3,T^{1,0}M,\theta)$ is pseudo-Einstein, then $P_{\theta}^\prime$ takes a simpler form:
$$P_\theta^\prime f = 4\Delta_b^2 f - 8 \text{Im}\left(\nabla^1(A_{11}\nabla^1 f)\right) - 4\text{Re}\left(\nabla^1(R\nabla_1 f)\right).$$
The computations of $P'_{\theta}\ln(\rho)$ shown in \cite{CCY2}, combined with the local expansion of the $\tau$ in Theorem \ref{texp} show that the Green's function $G_{\theta}$ of $A_{\theta}=\tau P'_{\theta}\tau$ has the following expansion:
$$G_{\theta}(x,y)=-\gamma_{3}\ln(d_{\theta}(x,y))+\mathcal{K}(x,y),$$
with $\mathcal{K}(x,y)\in L^{\infty}(M)$.\\

For the rest of the paper, $(M,T^{1,0}M,\theta)$ will always be assumed to be embeddable with $P'_{\theta}$ non-negative and $\ker P'_{\theta}=\R$.

\subsection{The Heisenberg group}
We identify the Heisenberg  group $\H$ with $\C\times \R\simeq\R^{3}$ with elements $w=(z,t)=(x+iy,t)\simeq(x,y,t)\in \R\times\R\times\R$ and group law
\begin{equation*}
w\cdot w'=(z,t)\cdot(z',t')=(z+z',t+t'+2\text{Im}(z\ol{z'}))\quad\forall\;  w,w'\in\H,
\end{equation*}
where $\text{Im}$ denotes the imaginary part of a complex number and $z\ol{z'}$ is the standard Hermitian inner product in $\C$. The dilations in $\H$ are
\begin{equation*}
\delta_{\lambda}:\H\to\H \qquad \delta_{\lambda}(z,t)=(\lambda z,\lambda^2 t)\quad \forall\;\lambda>0.
\end{equation*}
The natural distance that we will adopt in our setting is the Kor\'anyi distance, given by
\begin{equation*}
d_{\H}((z,t),(z',t'))=\left(|z-z'|^4+(t-t'-2\textnormal{Im}(z\ol{z'}))^2\right)^{\frac{1}{4}}
\end{equation*}
We denote by
\begin{equation*}
\Theta=\text{d}t+2\sum_{j=1}^{N}(x_i\text{d}y_i-y_i\text{d}x_i)
\end{equation*}
the standard contact form on $\H$ and by $dv_{\H}$ the volume form associated to $\Theta$. The Heisenberg group can be identified with the unit sphere in $\C^{2}$ minus a point through the Cayley transform $\Cay:\H\to S^{3}\setminus \{(0,0,0,-1)\}$ defined as follows
\begin{equation*}
\Cay(z,t)=\left(\frac{2z}{1+|z|^2+it},\frac{1-|z|^2-it}{1+|z|^2+it}\right).
\end{equation*}
On the unit sphere $S^{3}=\{\zeta\in\C^{2}:\; |\zeta|=1\}$ we consider the distance
\begin{equation*}
d_{S^{3}}(\zeta,\eta)^2=2|1-\zeta\ol{\eta}|,\quad \zeta,\eta\in \C^{2}
\end{equation*}
With this definition of $d_{S^{3}}$, the relation between the distance of two points $w=(z,t)$, $w'=(z',t')$ in $\H$ and the distance of their images $\Cay(w)$, $\Cay(w')$ in $S^{3}$, is given by
\begin{equation*}
d_{S^{3}}(\Cay(w),\Cay(w'))=d_{\H}(w,w')\left(\frac{4}{(1+|z|^2)^2+t^2}\right)^{\frac{1}{4}}\left(\frac{4}{(1+|z'|^2)^2+t'^2}\right)^{\frac{1}{4}}.
\end{equation*}
On $S^{3}$, we consider the standard contact form
\begin{equation*}
\theta_{0}=i\sum_{j=1}^{N+1}(\zeta_j\text{d}\ol{\zeta}_j-\ol{\zeta}_j\text{d}\zeta_j),
\end{equation*}
and we denote by $dv_{0}$ the volume form associated to $\theta_{0}$. With this notation we have that
$$(\Cay^{-1})^{*}\theta_{0}=|J_{\Cay}|^{\frac{1}{2}}\Theta,$$
where $|J_{\Cay}|=\frac{8}{[(1+|z|^{2})^{2}+t^{2}]^{2}}$ is the Jacobian of $\Cay$. For $h\in Aut(\H)$, we can parametrize their Jacobian $|J_{h}|$ as follows:
$$|J_{h}|=\frac{C}{||z|^{2}+it+2zw+\lambda|^{4}},$$
where $C>0$, $\lambda,w\in \C$ and $\text{Re}(\lambda)>|w|^{2}$. We also recall that $$Aut(S^{3})=\{k;k=\Cay\circ h \circ \Cay^{-1}, h\in Aut(\H)\}.$$
Hence,the Jacobian of $J_{k}$ can be parametrized as follow
$$|J_{k}|:=\frac{C}{|1-w\cdot \zeta|^{4}},$$
where $C>0$, $w\in \C^{2}$, $|w|<1$ and $\zeta\in S^{3}$. \\

We finish this section by this theorem regarding the pseudo-Hermitian normal coordinates:
\begin{theorem}[\cite{FS}]\label{thmfs}
Let $(M,\theta)$ be a pseudo-Hermitian manifold. Given $p\in M$, there exist neighborhoods $U$ of $p$ and $V$ of the origin of $\H$ and a diffeomorphism $\Psi:U\to V$ such that
\begin{itemize}
\item[i)] $(\Psi^{-1})^{*}\theta=(1+O_{1})\Theta.$
\item[ii)]$(\Psi^{-1})^{*}(\theta\wedge d\theta)=(1+O_{1})\Theta\wedge d\Theta$.
\end{itemize}
Where $O_{1}$ is a function satisfying $|O_{1}(x)|\leq C |x|$. 
\end{theorem}

\section{Properties of the mass and proof of Theorem \ref{thrm1}}
First, we start by determining, the change of the mass $m_{\theta}$ under a conformal change of the contact form $\theta\mapsto \theta_{F}=F^{\frac{1}{2}}\theta$.
\begin{proposition}\label{lm}
If $\theta_{F}=F^{\frac{1}{2}}\theta$, then 
$$m_{\theta_{F}}=m_{\theta}(x)+\frac{\gamma_{3}}{4}\ln(F(x))-\frac{2}{V_{F}}A^{-1}_{\theta}\tau F(x)+\frac{1}{V_{F}^{2}}\int_{M}FA^{-1}_{\theta}\tau F\ dv_{\theta}.$$
\end{proposition}
{\it Proof:}
We recall that based on our convention, the Green's function of the operator $A_{\theta}$ has the following properties:
$$\left\{\begin{array}{lll}
A_{\theta,x}G_{\theta}(x,y)=-\frac{1}{V}, \text{ for } x\not=y \\
G_{\theta}(x,y)+\gamma_{3}\ln(d_{\theta}(x,y)) \in L^{\infty}(M)\\
\int_{M}G_{\theta}(x,y)\ dv_{\theta}(y)=0.
\end{array}
\right.$$
We introduce then the function $$H_{\theta}(x,y,z):=G_{\theta}(x,y)-G_{\theta}(z,y).$$
Then one has:
$$A_{\theta,y}H_{\theta}=\delta_{x}-\delta_{z}.$$
Now, notice that by definition of the Green's function
$$A_{\theta,y}(H_{\theta_{F}}-H_{\theta})=0.$$
Moreover, $H_{\theta_{F}}-H_{\theta}$ is bounded. Hence, $H_{\theta_{F}}-H_{\theta}=constant$. Thus
$$G_{\theta_{F}}(x,y)-G_{\theta_{F}}(z,y)=G_{\theta}(x,y)-G_{\theta}(z,y)+C.$$
Integrating with respect to $dv_{\theta_{F}}(y)$ yields
$$A_{\theta}^{-1}\tau F(x)-A_{\theta}^{-1}\tau F(z)+CV_{F}=0.$$
Hence, $C=\frac{1}{V_{F}}A_{\theta}^{-1}\tau F(z)-A_{\theta}^{-1}\tau F(x))$. In particular
$$G_{\theta_{F}}(x,y)-G_{\tilde{\theta}}(z,y)=G_{\theta}(x,y)-G_{\theta}(z,y)+\frac{1}{V_{F}}A_{\theta}^{-1}\tau F(z)-A_{\theta}^{-1}\tau F(x)).$$
We integrate now with respect to $dv_{\theta_{F}}(z)$ to get
$$G_{\theta_{F}}(x,y)=G_{\theta}(x,y)-\frac{1}{V_{F}}A^{-1}_{\theta}\tau F(y)+\frac{1}{V_{F}^{2}}\int_{M}FA_{\theta}^{-1}\tau F \ dv_{\theta}-\frac{1}{V_{F}}A^{-1}_{\theta}\tau F(x).$$
Thus,
$$m_{\theta_{F}}(x)=m_{\theta}(x)+\frac{\gamma_{3}}{4}\ln(F(x))-\frac{2}{V_{F}}A^{-1}_{\theta}\tau F(x)+\frac{1}{V_{F}^{2}}\int_{M}FA^{-1}_{\theta}\tau F\ dv_{\theta}.$$
\hfill$\Box$

We also point out that a different proof of this result can be deduced from Lemma \ref{ap1} in the Appendix. A direct consequence of previous proposition is
$$\mathcal{M}_{\theta_{F}}=\int_{M}m_{\theta}F\ dv_{\theta}+\frac{\gamma_{3}}{4}\int_{M}F\ln(F)\ dv_{\theta}-\frac{1}{V_{F}}\int_{M}FA_{\theta}^{-1}\tau F\ dv_{\theta}.$$

\begin{corollary}[Theorem \ref{thrm1}]
Assume that $m_{\theta}$ is constant, then under the constraint $V_{F}=V$, one has:
$$\mathcal{M}_{\theta_{F}}-\mathcal{M}_{\theta}=\frac{\gamma_{3}}{4}\int_{M}F\ln(F)\ dv_{\theta}-\frac{1}{V}\int_{M}FA_{\theta}^{-1}\tau F\ dv_{\theta}.$$
In particular, on the standard sphere $(S^{3},\theta_{0})$,  one has
$$\mathcal{M}_{\theta_{F}}-\mathcal{M}_{\theta_{0}}=\frac{\gamma_{3}}{4}\int_{S^{3}}F\ln(F)\ dv_{\theta}-\frac{1}{V}\int_{S^{3}}FA_{\theta}^{-1}\tau F\ dv_{\theta}\geq 0.$$
with equality if an only if $F=|J_{k}|$ with $k\in Aut(S^{3})$ normalized to have volume $V$.
\end{corollary}
Notice that the last inequality follows from the LHLS inequality proved in \cite{Bran}. In fact, this  can be seen as the CR version of the spectral inequality in \cite{Mor1}.

\section{Proof of Theorem \ref{thrm2} part $( i)$}
We define the functional $J(\cdot,M):\L(M)\to \R$ by
$$J(F,M):=\int_{M}m_{\theta}F\ dv_{\theta}+\frac{\gamma_{3}}{4}\int_{M}F\ln(F)\ dv_{\theta}-\frac{1}{V_{F}}\int_{M}FA_{\theta}^{-1}\tau F\ dv_{\theta}.$$
In a similar way, for the Heisenberg group, we define the functional $J(\cdot,\H):\L_{c}(\H)\to \R$ by
$$J(f,\H):=\frac{\gamma_{3}}{4}\Big(\int_{\H}f\ln(f)dx-\frac{4}{V_{f}}\int_{\H}\int_{\H}f(x)\ln(\frac{1}{|xy^{-1}|})f(y)dxdy\Big),$$
and we let $$\mathcal{M}(\H,V)=\inf_{f\in \L_{c}(\H);V_{f}=V}J(f,\H).$$
We claim that 
\begin{equation}\label{m=m}
\mathcal{M}(\H,V)=\mathcal{M}([\theta_{0}],S^{3}).
\end{equation}
Indeed, from Theorem \ref{thrm1}, we have that $$\mathcal{M}([\theta_{0}],S^{3})= \mathcal{M}_{\theta_{0}}(S^{3}).$$
The equality in $(\ref{m=m})$, follows then from an easy computation starting from the LHLS inequality in $\H$, proved in $\cite{Bran,FL}$ and stated in the theorem below.
\begin{theorem}
For any measurable function $g:\H\to \R$ such that $g\geq 0$, $\int_{\H}g(x)\ dx=\omega_{3}:=2\pi^{2}$ and $\int_{\H}g\ln(1+|x|^{2})\ dx<\infty$, we have
$$\frac{2}{\omega_{3}^{2}}\int_{\H\times \H}\ln\Big(\frac{2}{|xy^{-1}|}\Big)g(x)g(y)\ dxdy\leq \frac{1}{\omega_{3}}\int_{\H}g\ln(g)\ dx+\ln(2).$$
with equality if and only if $g=(|J_{\Cay}|\circ h)|J_{h}|$ with $h\in Aut(\H)$.
\end{theorem}

Next, we calim that
\begin{equation}
\mathcal{M}([\theta],M)\leq \mathcal{M}(\H,V).
\end{equation}
But, in order to show this, we need an intermediate localization lemma. 

\begin{lemma}\label{coord}
Given $\varepsilon>0$, and $p\in M$, there exists $\delta>0$ and new coordinate system in $B_{\delta}(p)$ defined by a diffeomorphism $\phi$, such that
$$(\phi^{-1})^{*}(\theta\wedge d\theta)=\Theta\wedge d\Theta,$$
and
$$e^{-\varepsilon}\leq \frac{d_{\theta}(p,q)}{d_{\H}(\phi(p),\phi(q))}\leq e^{\varepsilon},$$
for all $q\in B_{\delta}(p)$.
\end{lemma}
{\it Proof:}
 First notice that using Theorem \ref{thmfs}, we have that for every $p\in M$, there exists $\delta>0$ and a diffeomorphism $\Psi:B_{\delta}(p)\to V$, where $V$ is a neighborhood of the origin in $\H$, such that
$$(\Psi^{-1})^{*}\theta=(1+O(\delta))\Theta,$$
and 
$$(\Psi^{-1})^{*}(\theta\wedge d\theta)=(1+O(\delta))\Theta\wedge d\Theta.$$
Now using Gray's theorem, we can find new coordinate systems in $\H$, defined by a diffeomerphism $\Phi$ such that
$$(\Phi^{-1})^{*}(1+O(\delta))(\Theta\wedge d\Theta)=\Theta \wedge d\Theta.$$
Since $(1+O(\delta))\Theta$ is close to $\Theta$ for $\delta$ small enough, and $\phi=\Phi\circ \Psi$, then given  $\varepsilon>0$, there exists $\delta>0$ such that
$$e^{-\varepsilon}\leq \frac{d_{\theta}(p,q)}{d_{\H}(\phi(p),\phi(q))}\leq e^{\varepsilon},$$
for all $q\in B_{\delta}(p)$.
\hfill$\Box$

\begin{lemma}\label{app}
Given $\varepsilon>0$, there exists $\delta>0$ such that, for $F\in\L(M)$ supported in $B_{\delta}(p)$, there exists $f\in 
\L_{c}(\H)$ compactly supported, such that 
$$|J(F,M)-J(f,\H)|\leq \varepsilon V_{F}.$$
Similarly, for any $f\in \L_{c}(\H)$, there exists $F\in \L(M)$ such that $F$ is supported in $B_{\delta}(p)$ and
$$|J(F,M)-J(f,H)|\leq \varepsilon V_{f}.$$
\end{lemma}
{\it Proof:}
Using the compactness of $M$ and a covering argument, we can always find $\delta>0$ such that for every $p\in M$, $B_{\delta}(p)$ is in a coordinate chart as described in Lemma \ref{coord}. So we fix $\varepsilon>0$. Taking $\delta>0$ even smaller if necessary, we can assume that
$$\Big|G_{\theta}(p,q)+\gamma_{3}\ln(|\phi(p)\phi(q)^{-1}|)-m_{\theta}(p)\Big|<\varepsilon.$$
Hence if $F$ is supported in $B_{\delta}(p)$, taking $U=\phi(B_{\delta}(p))$, $x=\phi(p)$, $y=\phi(q)$ and $f(x)=F(p)$, we have
\begin{align}
J(F,M)&=\frac{\gamma_{3}}{4}\int_{U}f\ln(f)\ dv_{\Theta}-\frac{\gamma_{3}}{V_{F}}\int_{U\times U}f(x)\ln(|xy^{-1}|)f(y)dv_{\Theta}(y) dv_{\Theta}(y)\notag\\
&\quad+\frac{1}{V_{F}}\int_{M\times M}F(p)\eta(p,q)F(q)\ dv_{\theta}(p)dv_{\theta}(q)\notag\\
&=J(f,\H)+\frac{1}{V_{F}}\int_{M\times M}F(p)\eta(p,q)F(q)\ dv_{\theta}(p)dv_{\theta}(q).
\end{align}
Hence,
$$|J(F,M)-J(f,\H)|\leq \varepsilon V_{F}.$$
In a similar way, the second assertion follows easily from the invariance of the functional $J(\cdot,\H)$ by the scaling 
\begin{equation}\label{scal}
f\mapsto \frac{1}{\lambda^{4}}f(\delta_{\frac{1}{\lambda}}\cdot),
\end{equation}
where $\delta_{\lambda}$ is the dilation in the Heisenberg group. So one can shrink the support and then lift it to a function on $M$ via the diffeomorphism $\phi$.
\hfill$\Box$
\begin{corollary}
Let $$\mathcal{M}_{\delta}([\theta],M):=\inf_{F\in \L(M); V_{F}=V \text{ and } supp(F)\subset B_{\delta}(p)}J(F,M).$$ 
Then, one has
$$\lim_{\delta\to 0}\mathcal{M}_{\delta}([\theta],M)=\mathcal{M}(\H,V).$$
In particular, $$\mathcal{M}([\theta],M)\leq \mathcal{M}(\H,V).$$
\end{corollary}

\section{Proof of Theorem \ref{thrm2} part $(ii)$}

\subsection{Concentration and Improved LHLS inequality}

We move now to proving some regularity estimates that follow from the structure of the Green's function of $A_{\theta}$.
\begin{proposition}
For any $\varepsilon>0$, there exists $C_{\varepsilon}>0$ such that for all $F\in \L(M)$, 
\begin{equation}\label{est1}
\|A^{-1}_{\theta}\tau F\|_{\infty}\leq (1+\varepsilon)\frac{\gamma_{3}}{4}\int_{M}F\ln(F)\ dv_{\theta}+C_{\varepsilon}\Big(\int_{M}F\ dv_{\theta}+1\Big).
\end{equation}
\end{proposition}
{\it Proof:}
Recall that $$G_{\theta}=-\gamma_{3}\ln(d_{\theta}(x,y))+\mathcal{K}(x,y),$$
where $\mathcal{K}$ is a bounded smooth kernel. Thus for $r>0$ and small,
\begin{align}
A^{-1}_{\theta}F(x)&=-\gamma_{3}\int_{M}\ln(d_{\theta}(x,y))F(y)\ dv_{\theta}(y)+\int_{M}\mathcal{K}(x,y)F(y)\ dv_{\theta}(y)\notag\\
&=-\gamma_{3}\int_{B_{r}(x)}\ln(|xy^{-1}|)F(y)dy+\int_{M}\tilde{K}(x,y)F(y)\ dv_{\theta}(y).\notag
\end{align}
The second term is naturally bounded by $C\int_{M}F\ dv_{\theta}$. For the first term, we have for $0<\delta<4$
\begin{align}
-4\int_{B_{r}(x)}\ln(|xy^{-1}|)F(y)dy&=-4\int_{\ln(F)\geq-(4-\delta)\ln(|xy^{-1}|)}\ln(|xy^{-1}|)F(y)\ dy\notag\\
&\quad-4\int_{\ln(F)\leq-(4-\delta)\ln(|xy^{-1}|)}\ln(|xy^{-1}|)F(y)\ dy\notag\\
&\leq \frac{4}{4-\delta}\int_{M}F\ln(F)\ dv_{\theta}-4\int_{B_{r}(x)}\frac{1}{|xy^{-1}|^{4-\delta}}\ln(|xy^{-1}|)\ dy
\end{align}
The second term in the last inequality is clearly bounded for all $0<\delta<4$. Hence, $(\ref{est1})$  is proved by choosing $\delta$ arbitrarily close to $4$.
\hfill$\Box$

Next we state this useful fact that follows directly from the Cauchy-Schwartz inequality.
\begin{lemma}\label{cs}
We let $\sigma_{F}=\frac{V_{F}}{V}$, for $F\in \L(M)$ and assume that $Q$ and $R$ satisfy the inequality
$$\frac{1}{V}\int_{M}FA^{-1}_{\theta}\tau F \ dv_{\theta}\leq \frac{\gamma_{3}}{4}\int_{M}F\ln(F)\ dv_{\theta}+C, \text{ if } V_{F}=V.$$
Then,
$$\frac{1}{V}\int_{M}QA^{-1}_{\theta}\tau R\ dv_{\theta}\leq \frac{\gamma_{3}}{8}\Big(\sigma_{R}\int_{M}Q\ln\Big(\frac{Q}{\sigma_{Q}}\Big)\ dv_{\theta}+\sigma_{Q}\int_{M}R\ln\Big(\frac{R}{\sigma_{R}}\Big)\ dv_{\theta}\Big)+C\sigma_{Q}\sigma_{R}.$$
\end{lemma}
Now, we have the necessary tools to prove a weak logarithmic HLS inequality:
\begin{proposition}\label{weak}
There exists $C>0$, depending on $M$ and $V$, such that if $V_{F}=V$ one has
\begin{equation}\label{hardyg}
\frac{1}{V}\int_{M}FA^{-1}_{\theta}\tau F \ dv_{\theta}\leq \frac{\gamma_{3}}{4}\int_{M}F\ln(F)\ dv_{\theta}+C.
\end{equation}

\end{proposition}
{\it Proof:}
First, notice that if we apply Lemma \ref{app} with $\varepsilon=1$, then one has that there exists $\delta>0$ such that for $F$ supported in $B_{\delta}(p)$, we have
$$J(F,M)\geq \mathcal{M}(\H,V)-V_{F}.$$
Hence, by the boundedness of $m_{\theta}$, we have that $(\ref{hardyg})$ holds for functions compactly supported in $B_{\delta}(p)$. \\
Next, we consider a covering of $M$ by closed sets $(U_{i})_{1\leq i \leq N}$ such that if $U_{i}\cap U_{j}=\emptyset$ then $d_{\theta}(U_{i},U_{j})>\varepsilon$. We can also assume that the sets $U_{i}$ are small, in a way that if $U_{i}\cap U_{j}\not= \emptyset$ then there exists a ball $B_{\delta}(p)$ containing them both. We set $F_{i}=\chi_{U_{i}}F$. So one can use the fact that the Green's function is regular away from the diagonal, in order to write
\begin{align}
\int_{M}FA^{-1}_{\theta}\tau F\ dv_{\theta}&=\sum_{U_{i}\cap U_{j}=\emptyset}\int_{M}F_{i}A^{-1}_{\theta}\tau F_{j}\ dv_{\theta}+\sum_{U_{i}\cap U_{j}\not=\emptyset}\int_{M}F_{i}A^{-1}_{\theta}\tau F_{j}\ dv_{\theta}\notag\\
&\leq V_{F}^{2}C_{1}+\sum_{U_{i}\cap U_{j}\not=\emptyset}\int_{M}F_{i}A^{-1}_{\theta}\tau F_{j}\ dv_{\theta}\notag\\
&=V_{F}^{2}C_{1}+II
\end{align}
Therefore, one should focus on the term $II$. But using Lemma \ref{cs}, one has
$$
II\leq \sum_{U_{i}\cap U_{j}\not=\emptyset}\Big(\frac{V_{F_{i}}\gamma_{3}}{8}\int_{M}F_{i}\ln\Big(\frac{F_{i}}{\sigma_{F_{i}}}\Big)\ dv_{\theta}+\frac{V_{F_{j}}\gamma_{3}}{8}\int_{M}F_{j}\ln\Big(\frac{F_{j}}{\sigma_{F_{j}}}\Big)\ dv_{\theta}+C\frac{V_{F_{i}}V_{F_{j}}}{V}\Big).
$$
Then, applying Jensen's inequality for the function $t\mapsto t\ln(t)$, under the assumption that $V_{F}=V$, yields
$$\int_{M}F_{i}\ln\Big(\frac{F_{i}}{\sigma_{j}}\Big)\ dv_{\theta}\geq 0.$$
Thus,
\begin{align}
&\sum_{U_{i}\cap U_{j}\not=\emptyset}\Big(\frac{V_{F_{i}}\gamma_{3}}{8}\int_{M}F_{i}\ln\Big(\frac{F_{i}}{\sigma_{F_{i}}}\Big)\ dv_{\theta}+\frac{V_{F_{j}}\gamma_{3}}{8}\int_{M}F_{j}\ln\Big(\frac{F_{j}}{\sigma_{F_{j}}}\Big)\ dv_{\theta}+C\frac{V_{F_{i}}V_{F_{j}}}{V}\Big)\label{volcomp}\\
&\leq \sum_{i,j}\Big(\frac{V_{F_{i}}\gamma_{3}}{8}\int_{M}F_{i}\ln\Big(\frac{F_{i}}{\sigma_{F_{i}}}\Big)\ dv_{\theta}+\frac{V_{F_{j}}\gamma_{3}}{8}\int_{M}F_{j}\ln\Big(\frac{F_{j}}{\sigma_{F_{j}}}\Big)\ dv_{\theta}+C\frac{V_{F_{i}}V_{F_{j}}}{V}\Big)\notag\\
&=\frac{\gamma_{3}}{4}V\int_{M}F\ln(F)\ dv_{\theta}-\frac{\gamma_{3}}{4}V\sum_{j}V_{F_{j}}\ln(V_{F_{j}})+C_{2}.\notag
\end{align}
Taking $C_{3}=-\min t\ln(t)$, yields
$$\int_{M}FA^{-1}_{\theta}\tau F\ dv_{\theta}\leq \frac{\gamma_{3}}{4}V\int_{M}F\ln(F)\ dv_{\theta}+\frac{\gamma_{3}}{4}C_{3}NV+C_{2}+C_{1}V^{2},$$
which is the desired inequality.
\hfill$\Box$\\

As it was noted in the proof of Lemma \ref{app}, the functional $J(\cdot,\H)$ is invariant under the scaling $(\ref{scal})$, which leaves the volume, or the $L^{1}$-norm invariant. This hints to a concentration phenomena that can happen locally for the functional $J(\cdot,M)$. So we start investigating the effect of concentration on the functional $J(\cdot,M)$.
\begin{definition}
We say that a sequence $(F_{j})_{j\in \N} \in L^{1}(M)$ is a concentrating sequence, if there exists a sequence of points $p_{j}\in M$ and numbers $\delta_{j}\to 0$ such that
$$\int_{B_{\delta_{j}}(p_{j})}F_{j}\ dv_{\theta}\geq (1-\delta_{j})V_{F_{j}}.$$
\end{definition}
We then have this lower bound on the energy of concentrating sequences:
\begin{proposition}\label{con}
Let $(F_{j})_{j\in \N}$ be a concentrating sequence in $\L(M)$ with constant volume $V$. Then $$\liminf _{j\to \infty}J(F_{j},M)\geq \mathcal{M}(\H,V).$$
\end{proposition}
{\it Proof:}
Let $(F_{j})_{j\in \N}$ such a sequence and consider $\chi_{j}$ the characteristic function of $B_{\delta_{j}}(p_{j})$. We define then $Q_{j}=\chi_{j}F_{j}$ and $R_{j}=F_{j}-Q_{j}$. One then has from Lemma \ref{app}, that for $j$ large enough,
$$J(Q_{j},M)\geq \mathcal{ M}(\H,V_{Q_{j}})-\varepsilon V_{Q_{j}}.$$
Using Lemma \ref{cs}, we have
$$\int_{M}Q_{j}A^{-1}_{\theta}\tau R_{j}\ dv_{\theta}\leq\frac{\gamma_{3}}{8}\Big(V_{R_{j}}\int_{M}Q_{j}\ln\Big(\frac{Q_{j}}{\sigma_{Q_{j}}}\Big)\ dv_{\theta}+V_{Q_{j}}\int_{M}R_{j}\ln\Big(\frac{R_{j}}{\sigma_{R_{j}}}\Big)\ dv_{\theta}\Big)+C\sigma_{Q_{j}}\sigma_{R_{j}}.$$
Since $V=V_{F_{j}}$, one has
\begin{align}
VJ(F_{j},M)&=(V_{Q_{j}}+V_{R_{j}})\int_{M}m_{\theta}(Q_{j}+R_{j})-\int_{M}(Q_{j}+R_{j})A^{-1}_{\theta} \tau(Q_{j}+R_{j})\ dv_{\theta}\notag\\
&\quad+(V_{Q_{j}}+V_{R_{j}})\frac{\gamma_{3}}{4}\Big(\int_{M}Q_{j}\ln(Q_{j})\ dv_{\theta}+\int_{M}R_{j}\ln(R_{j})\ dv_{\theta}\Big)\notag\\
&=V_{Q_{j}}J(Q_{j},M)+\frac{\gamma_{3}}{4}V_{R_{j}}\int_{M}R_{j}\ln(R_{j})\ dv_{\theta}-\int_{M}R_{j}A^{-1}_{\theta}\tau R_{j}\ dv_{\theta}\notag\\
&\quad+V_{Q_{j}}\frac{\gamma_{3}}{4}\int_{M}R_{j}\ln(R_{j})\ dv_{\theta}+V_{R_{j}}\frac{\gamma_{3}}{4}\int_{M}Q_{j}\ln(Q_{j})\ dv_{\theta}-2\int_{M}Q_{j}A^{-1}_{\theta}\tau R_{j}\ dv_{\theta} \notag\\
&\quad+\int_{M}(V_{Q_{j}}R_{j}+V_{R_{j}}F_{j})m_{\theta}\ dv_{\theta}.\notag
\end{align}
Therefore,
\begin{align}
VJ(F_{j},M)&\geq V_{Q_{j}}J(Q_{j},M)+\frac{\gamma_{3}}{4}\Big(V_{R_{j}}^{2}\ln(\sigma_{R_{j}})+V_{R_{j}}V_{Q_{j}}\Big(\ln(\sigma_{R_{j}})+\ln(\sigma_{Q_{j}})\Big)\Big)\notag\\
&\quad-\frac{CV_{R_{j}}^{2}}{V}-\frac{2CV_{R_{j}}V_{Q_{j}}}{V}-V_{R_{j}}(V_{Q_{j}}+V)\max|m_{\theta}|\notag\\
&\geq V_{Q_{j}}J(Q_{j},M)+\frac{\gamma_{3}}{4}V_{R_{j}}V\ln(V_{R_{j}})-C'V_{R_{j}}.\notag
\end{align}
Thus,
$$J(F_{j},M)\geq \frac{V_{Q_{j}}}{V}\Big(\mathcal{M}(\H,V_{Q_{j}})-\varepsilon_{j}V_{Q_{j}}\Big)+\frac{\gamma_{3}}{4}V_{R_{j}}\ln(V_{R_{j}})-C'V_{R_{j}}.$$
So, passing to the limit, yields the desired result.
\hfill$\Box$\\

Since concentration tends to localize the problem in a way that it becomes similar to the Heisenberg case, one expects to obtain an improved logarithmic HLS inequality in the case of absence of concentration and this can be quantified by the following:
\begin{lemma}\label{imp}
Fix $0<\delta<1$, then there exists $C(\delta,M,\theta)>0$ such that for any $F\in \L(M)$ satisfying $V_{F}=V$ and 
$$\int_{B_{\delta}(x)}F\ dv_{\theta}<(1-\delta)V,$$
for all $x\in M$, we have
$$ (1-\delta)\frac{\gamma_{3}}{4}\int_{M}F\ln(F)\ dv_{\theta}+C \geq \frac{1}{V}\int_{M}FA^{-1}_{\theta}\tau F\ dv_{\theta}.$$
\end{lemma}
{\it Proof:}
The proof is similar to the one of Proposition \ref{weak}. So we consider a covering of $M$ by closed sets $(U_{i})_{1\leq i \leq N}$ such that if $U_{i}\cap U_{j}=\emptyset$ then there exists $\varepsilon>0$ such that $d_{\theta}(U_{i},U_{j})>\varepsilon$. We can also assume that the sets $U_{i}$ are small, in a way that if $U_{i}\cap U_{j}\not= \emptyset$ then there exist a number $0<\delta'<\delta$ and a ball $B_{\delta'}(p)$ containing them both. We add the extra condition that for all $1\leq j\leq N$, there exists $x_{j}\in M$  such that $U_{j}\subset B_{\frac{\delta}{2}}(x_{j})$. With the notations as above, we notice now that
$$\sum_{i;U_{j}\cap U_{i}\not=\emptyset}V_{F_{i}}\leq \int_{B_{\delta}(x_{j})}F\ dv_{\theta}\leq (1-\delta)V_{F},$$
Using (\ref{volcomp}), we get 
\begin{align}
\int_{M}FA_{\theta}^{-1}\tau F \ dv_{\theta}&\leq \frac{\gamma_{3}(1-\delta)V_{F}}{4}\sum_{i}\int_{M}F_{i}\ln\Big(\frac{F_{i}}{\sigma_{F_{i}}}\Big)\ dv_{\theta}+C+C_{1}V^{2}\notag\\
&\leq \frac{\gamma_{3}(1-\delta)V_{F}}{4}\int_{M}F\ln(F)\ dv_{\theta}+\tilde{C}.\notag
\end{align}
\hfill$\Box$

\subsection{Sub-critical Approximation}
Consider the modified functional
$$J_{\varepsilon}(F,M)=\int_{M}m_{\theta}F\ dv_{\theta}+\frac{\gamma_{3}}{4}\int_{M}F\ln(F)\ dv_{\theta}-\frac{(1-\varepsilon)\lambda_{1}^{\varepsilon}}{V_{F}}\int_{M}FA_{\theta}^{-1-\varepsilon}\tau F\ dv_{\theta},$$
where $\lambda_{1}$ is the first non-zero eigenvalue of $A_{\theta}$ and $\varepsilon>0$.
\begin{lemma}
There exists $F_{\varepsilon}\in C^{\infty}(M)$ that minimizes the functional $J_{\varepsilon}$. That is,
$$\inf_{f\in \L(M); V_{F}=V}J_{\varepsilon}(F,M)=J_{\varepsilon}(F_{\varepsilon},M).$$
\end{lemma}
{\it Proof:}
First notice that
$$J_{\varepsilon}(F,M)\geq \int_{M}m_{\theta}F\ dv_{\theta}+\frac{\gamma_{3}}{4}\int_{M}F\ln(F)\ dv_{\theta}-\frac{(1-\varepsilon)}{V_{F}}\int_{M}FA_{\theta}^{-1}\tau F\ dv_{\theta}.$$
Now using Proposition \ref{weak}, we have
\begin{align}
J_{\varepsilon}(F,M)&\geq \frac{\varepsilon\gamma_{3}}{4}\int_{M}F\ln(F)\ dv_{\theta}-C(1-\varepsilon)+\inf_{M}(m_{\theta})V\notag\\
&\geq \frac{\varepsilon\gamma_{3}}{4}\int_{M}F\ln(F)\ dv_{\theta}-C_{1}.
\end{align}
Therefore, if $(F_{k})_{k\in\N}$ is a minimizing sequence for $J_{\varepsilon}$, then $\int_{M}F_{k}\ln(F_{k})\ dv_{\theta}$ is bounded above independently of $k$.
To finish our argument, we recall the following useful result \cite[Lemma 2.11]{Ok1}.
\begin{lemma}
We fix a continuous convex function $G:[0,\infty)\to \R$ such that 
$$\lim_{t\to \infty}\frac{G(t)}{t}=\infty.$$
Consider a sequence $(F_{k})_{k\in \N}$ of non-negative measurable functions in $M$ such that
$$\sup_{k\in \N}\int_{M}G(F_{k})dv=S<\infty.$$
Then after passing to a subsequence, there exists $F\in L^{1}(M)$ such that $F_{k}\to F$ weakly and
$$\int_{M}G(F)dv\leq \liminf_{k\to \infty}G(F_{k})dv.$$
\end{lemma}
Using this Lemma for $G(t)=t\ln(t)$ and the sequence $(F_{k})_{k\in\N}$, we have the existence of $F_{\varepsilon}\in \L(M)$ such that
$F_{k}\to F_{\varepsilon}$ weakly in $L^{1}(M)$ and $$\int_{M}F_{\varepsilon}\ln(F_{\varepsilon})\ dv_{\theta} \leq \liminf_{k\to \infty} \int_{M}F_{k}\ln(F_{k})\ dv_{\theta}.$$
Since $F_{k}$ is bounded in $\L(M)$, we have from Proposition \ref{weak}, that $A^{-1}_{\theta}\tau F_{k}$ is uniformly bounded in $L^{p}(M)\cap \P$ for all $1\leq p \leq \infty$. Now by ellipticity of $A_{\theta}$ on $\P$, we have that $A_{\theta}^{-\varepsilon}\tau$ is a pseudo-differential operator of order $-4\varepsilon$. Hence, $A^{-1-\varepsilon}_{\theta}\tau F_{k}$ is uniformly bounded in $W^{p,4\varepsilon}(M)$. Taking $p>\frac{4}{3\varepsilon}$, we see that $(A^{-1-\varepsilon}_{\theta}\tau F_{k})_{k\in \N}$ is compact in $C(M)$. Therefore, we can extract a convergent subsequence, that we still denote by $(F_{k})_{k\in \N}$ such that $A^{-1-\varepsilon}_{\theta}\tau F_{k}\to A^{-1-\varepsilon}_{\theta}\tau F_{\varepsilon}$, since  $F_{k}\to F_{\varepsilon}$ weakly in $L^{1}(M)$. Hence,
$$\int_{M}F_{k}A^{-1-\varepsilon}_{\theta}\tau F_{k}\ dv_{\theta}\to \int_{M}F_{\varepsilon}A^{-1-\varepsilon}_{\theta}\tau F_{\varepsilon}\ dv_{\theta}.$$ 
Therefore, 
$$\lim_{k\to \infty}J_{\varepsilon}(F_{k},M)\geq J_{\varepsilon}(F_{\varepsilon},M).$$
Let us show that $F_{\varepsilon}$ is bounded below by a positive constant. Consider a bounded function $H$ such that $F_{\varepsilon}+H\in \L(M)$ and $\int_{M}H\ dv_{\theta}=0$, then we have
\begin{align}
\mathcal{D}&:=J_{\varepsilon}(F_{\varepsilon}+H,M)-J_{\varepsilon}(F_{\varepsilon},M)\notag\\
&=\int_{M}m_{\theta}H\ dv_{\theta}+\frac{\gamma_{3}}{4}\int_{M}(F_{\varepsilon}+H)\ln(F_{\varepsilon}+H)-F_{\varepsilon}\ln(F_{\varepsilon})\ dv_{\theta}\notag\\
&\quad-\frac{2(1-\varepsilon)\lambda_{1}^{\varepsilon}}{V}\int_{M}HA^{-1-\varepsilon}_{\theta}\tau F_{\varepsilon}\ dv_{\theta}-\frac{2(1-\varepsilon)\lambda_{1}^{\varepsilon}}{V}\int_{M}HA^{-1-\varepsilon}_{\theta}\tau H\ dv_{\theta}.
\end{align}
Hence,
$$\mathcal{D}\leq C\int_{M}|H|\ dv_{\theta}+\frac{\gamma_{3}}{4}\int_{M}(F_{\varepsilon}+H)\ln(F_{\varepsilon}+H)-F_{\varepsilon}\ln(F_{\varepsilon})\ dv_{\theta},$$
where $C=\|m_{\theta}\|_{\infty}+\frac{(1-\varepsilon)\lambda_{1}^{\varepsilon}}{V}\|A^{-1-\varepsilon}_{\theta}\tau F_{\varepsilon}\|_{\infty}$. Notice that from the mean value theorem, one has 
$$(t+s)\ln(t+s)-t\ln(t)<s(1+\ln(t+s))<0,$$
whenever $t,s>0$ and $t+s<e^{-1}$ or $t>0$, $s<0$ and $t+s>e^{-1}$. So we consider the two sets $W:=\{x;F_{\varepsilon}>e^{-1}\}$ and $W_{N}:=\{ x;F_{\varepsilon}(x)<e^{-N}\}$. Notice that by the mean value theorem, we have that $W$ has positive measure. So we assume for the sake of contradiction that also $W_{N}$ has positive measure and construct the function $H$ such that $\int_{M}H\ dv_{\theta}=0$ with
$$\left\{\begin{array}{lll}
0<H<F-e^{-1}& \text{ on } W\\
F-e^{-N}<H<0& \text{ on } W_{N}\\
H=0 &\text{ on } M\setminus (W\cup W_{N}).
\end{array}
\right.
$$
Notice that in this case we have 
\begin{align}
\mathcal{D}&\leq -N\frac{\gamma_{3}}{4}\int_{W_{N}}H\ dv_{\theta}+C\int_{M}|H|\ dv_{\theta}\notag\\
&\leq (C-N\frac{\gamma_{3}}{8})\int_{M}|H|\ dv_{\theta},
\end{align}
which yields a contradiction for $N$ big enough. Hence, $F_{\varepsilon}$ is bounded from below. Now the Euler-Lagrange equation for the constraint minimization of $J_{\varepsilon}$ yields the equation
\begin{equation}
m_{\theta}+\frac{\gamma_{3}}{4}(\ln(F_{\varepsilon})+1)-\frac{2(1-\varepsilon)\lambda_{1}^{\varepsilon}}{V}A^{-1-\varepsilon}_{\theta}\tau F_{\varepsilon}=\lambda_{\varepsilon},
\end{equation}
where $\lambda_{\varepsilon}$ is the constant coming from the Lagrange multiplier. Therefore, by ellipticity of $A_{\theta}$ restricted to $\P$ and smoothness of $m_{\theta}$, we get the smoothness of $F_{\varepsilon}$.
\hfill$\Box$\\

At this stage, we have the required ingredients to finish the proof of Theorem \ref{thrm2}. The idea is to extract a convergent subsequence of $F_{\varepsilon}$ when $\varepsilon \to 0$. Notice that we have $$\lim_{\varepsilon\to 0}J(F_{\varepsilon},M)=\mathcal{M}([\theta],M).$$
So $(F_{\varepsilon})_{\varepsilon>0}$ is a minimizing sequence that we need to show its convergence. But since $\mathcal{M}([\theta],M)< \mathcal{M}([\theta_{0}],S^{3})$, it follows from Proposition \ref{con} that $(F_{\varepsilon})_{\varepsilon>0}$ does not concentrate. We combine then Lemma \ref{imp} and the boundedness of $J(F_{\varepsilon},M)$ to get
$$\int_{M}F_{\varepsilon}\ln(F_{\varepsilon})\ dv_{\theta}\leq C,$$
with $C$ is independent of $\varepsilon$. This, combined with
$$J_{\varepsilon}(1,M)\geq J_{\varepsilon}(F_{\varepsilon},M)\geq \mathcal{M}([\theta],M),$$
yields the uniform boundednes of $\int_{M}F_{\varepsilon}A^{-1-\varepsilon}_{\theta}\tau F_{\varepsilon}\ dv_{\theta}.$

We recall now that $F_{\varepsilon}$ satisfies the equation
\begin{equation}\label{lag}
m_{\theta}+\frac{\gamma_{3}}{4}\Big(\ln(F_{\varepsilon})+1\Big)-\frac{2(1-\varepsilon)\lambda_{1}^{\varepsilon}}{V}A^{-1-\varepsilon}_{\theta}\tau F_{\varepsilon}=\lambda_{\varepsilon}.
\end{equation}
The Lagrange multiplier $\lambda_{\varepsilon}$ can be obtained by multiplying $(\ref{lag})$ by $F_{\varepsilon}$ and then integrating:
$$
\lambda_{\varepsilon}V=\int_{M}m_{\theta}F_{\varepsilon}\ dv_{\theta}+\frac{\gamma_{3}}{4}\int_{M}F_{\varepsilon}\ln(F_{\varepsilon})\ dv_{\theta}-\frac{2(1-\varepsilon)\lambda_{1}^{\varepsilon}}{V}\int_{M}F_{\varepsilon}A^{-1-\varepsilon}_{\theta}\tau F_{\varepsilon}\ dv_{\theta}+\frac{\gamma_{3}}{4}V.
$$
Hence, $\lambda_{\varepsilon}$ is uniformly bounded. Using Proposition \ref{weak}, we have that $(A_{\theta}^{-1}F_{\varepsilon})_{\varepsilon}$ is uniformly bounded in $C(M)\cap \P$. A boot-strap argument for equation $(\ref{lag})$ provides us with the smoothness of $F_{\varepsilon}$. So if we set $u_{\varepsilon}=A^{-1}_{\theta}\tau F_{\varepsilon}$, then one has
$$\int_{M}u_{\varepsilon}A_{\theta}u_{\varepsilon} \ dv_{\theta}\leq C.$$
Therefore, from the Moser-Trudinger inequality in \cite{CaYa1}, we have that $u_{\varepsilon}$ is uniformly bounded in $W^{2,2}(M)$ and $e^{u_{\varepsilon}}$ is uniformly bounded in $L^{p}(M)$ for all $1\leq p <\infty$. Since the family $(A^{-\varepsilon}_{\theta})_{\varepsilon}$ is uniformly bounded in $W^{2,2}(M)\cap \P$, we have the uniform boundedness in $W^{2,2}(M)\cap \P$ of $v_{\varepsilon}:=A^{-\varepsilon}_{\theta}u_{\varepsilon}$. Again, using the Moser-Trudinger inequality, we get that $e^{v_{\varepsilon}}$ is uniformly bounded in $L^{p}(M)$. But since
$$F_{\varepsilon}=e^{R_{\varepsilon}} e^{\frac{8(1-\varepsilon)\lambda_{1}^{\varepsilon}}{\gamma_{3}V}v_{\varepsilon}},$$
and $R_{\varepsilon}$ is uniformly bounded in $L^{\infty}(M)$, we have the uniform boundeddness of $F_{\varepsilon}$ in $L^{2}(M)$. So, using the regularizing effect of $A^{-1}_{\theta}$, we see that $(A_{\theta}^{-1-\varepsilon}\tau F_{\varepsilon})_{\varepsilon}$ is compact in $C(M)\cap \P$.  Therefore, we can extract a convergent subsequence of $(F_{\varepsilon})_{\varepsilon}$ that we denote by $(F_{\varepsilon_{k}})_{k\in \N}$ such that $F_{\varepsilon_{k}}\to F_{0}$ in $C(M)$ and via a diagonal process, we get that
$$J(F_{0},M)=\inf_{F\in \L(M);V_{F}=V}J(F,M).$$

\hfill$\Box$\\

It is then easy to see that if $\mathcal{M}([\theta],M)$ is achieved by a function $F_{0}$, then it satisfies the equation
$$m_{\theta}+\frac{\gamma_{3}}{4}\Big(\ln(F_{0})+1\Big)-\frac{2}{V}A_{\theta}^{-1}\tau F_{0}=\lambda_{0}.$$
In particular, a boot-strapping argument yields the regularity of $F_{0}$ and if $m_{\theta}\in \P$ then so is $\ln(F_{0})$. Using Proposition \ref{lm}, one sees that $m_{\theta_{F}}$ is constant. In fact, we have
$$m_{\theta_{F_{0}}}=\frac{\mathcal{M}([\theta],M)}{V}.$$

\section*{Appendix}
\appendix
\section{The regularized Zeta function and the Mass}
In this section we will establish a link between the total mass and the regularized Zeta function. We want to point out that in the Riemannian case, this link was established in $\cite[Section 5]{Mor2}$ for a general pseudo-differential operator having a leading term of $\Delta_{g}^{\frac{d}{2}}$ where $d$ is the dimension of the manifold without the mention of concept of mass. In \cite{Ok1}, the author introduced the mass. Both these two proves rely on the heat kernel estimates and an explicit expansion of the Green's function of the fractional power of operator. In our case, we avoid the use of the heat kernel expansion since the operator $A_{\theta}$ does not have an invertible symbol as an operator in $\Psi_{H}(M)$ ( but $A_{\theta}$ is a generalized Toeplitz operator which is invertible in the Toeplitz algebra \cite{Bout}). Therefore one cannot use the Volterra calculus and the heat kernel expansion developed in \cite{Pon1}. Our proof relies on the non-commutative residue introduced in \cite{Pon}.\\

We will be using the same notations as \cite{Pon1,Pon}. We consider the operator $P_{0}\in \Psi_{H}(M)$ with Schwartz kernel $K_{0}\in \mathcal{K}^{0}(M\times M)$ such that in a coordinate patch around $x\in M$ we have $K_{0}(x,y)=-\gamma_{3}\ln(d_{\theta}(x,y))$ (Here we disregarded the factor related to the Jacobian of the change of coordinates for the sake of notation). Then we have
$$\mathcal{M}_{\theta}=\int_{M}\lim_{x\to y}(G_{\theta}(x,y)-K_{0}(x,y))dv_{\theta}(x)$$
Therefore, since $A_{\theta}^{-1}\tau-P_{0}$ is a trace class operator,
$$\mathcal{M}_{\theta}=TR(A_{\theta}^{-1}\tau-P_{0}).$$
We consider now the holomorphic family $s\mapsto A_{\theta}^{-s}$ defined in a neighborhood of zero, where for $\Re(s)>0$ we have
$$A_{\theta}^{-s}:=q(s)\int_{0}^{\infty}t^{-s}(A_{\theta}+r+t)^{-1}dt$$
where $q(s)=\frac{1}{\int_{0}^{\infty}t^{-s}(1+t)^{-1}dt}$ and $r:\P\to \ker A_{\theta}$ is the $L^{2}$-orthogonal projection. Notice that $A_{\theta}^{-s}$ is defined on $\P$ and can be extended by 0 to $\mathcal{P}^{\perp}$. We also have $ord(A_{\theta}^{-s})=-4s$.\\
Notice that with the previous notation $\lim_{s\to 0^{+}}A_{\theta}^{-s}u=u-r(u)$ for all $u\in \P$. So we let $T_{s}=A_{\theta}^{-s}\oplus \tau^{\perp}\oplus r$. We will be using this family as a gauging for $A_{\theta}^{-1}\tau$ since for $\Re(s)>0$ and small, $T_{s}A_{\theta}^{-1}\tau$ is a trace class operator.\\

Recall that, \cite[Proposition 3.17]{Pon}, $TR(T_{s}A_{\theta}^{-1}\tau)$ has a simple pole at $s=0$ and the residue at this pole is
$$Res_{s=0}\Big(TR(T_{s}A_{\theta}^{-1}\tau)\Big)=-Res(A^{-1}_{\theta}\tau)=-\gamma_{3}V.$$
Similarly $TR(T_{s}P_{0})$ has the same residue at the pole $s=0$. Hence,
\begin{align}
\lim_{s\to 0}TR(T_{s}A_{\theta}^{-1}\tau)-TR(T_{s}P_{0})&=\lim_{s\to 0}TR(T_{s}A_{\theta}^{-1}\tau)-\frac{(-\gamma_{3}V)}{s}+\frac{(-\gamma_{3}V)}{s}-TR(T_{s}P_{0})\notag\\
&=\tilde{Z}(1)+c,\label{las}
\end{align}
where $$c=\lim_{s\to 0}\frac{(-\gamma_{3}V)}{s}-TR(T_{s}P_{0}),$$
is a constant that might depend on $V$. The last equality in $(\ref{las})$ follows from the fact that $T_{s}A^{-1}_{\theta}\tau=A_{\theta}^{-s}A_{\theta}^{-1}\tau$. But
$$\lim_{s\to 0}TR(T_{s}A_{\theta}^{-1}\tau-T_{s}P_{0})=TR(A_{\theta}^{-1}\tau-P_{0}))=\mathcal{M}_{\theta}.$$
Therefore, 
$$\tilde{Z}(1)=\mathcal{M}_{\theta}-c.$$

\section{Geometric CR mass}
In this section we will add a geometric correction to the mass that makes it independent of the point on the sphere. First, we start by the following:
\begin{lemma}\label{ap1}
$$A_{\tilde{\theta}}^{-1}\tau_{F}f=A_{\theta}^{-1}\tau(Ff)-\frac{\int_{M}Ff\ dv_{\theta}}{V_{F}}A_{\theta}^{-1}\tau(F)-a_{1}+a_{2}$$
where $a_{1}=\frac{\int_{M}FA_{\theta}^{-1}\tau(Ff)\ dv_{\theta}}{V_F}$  $a_{2}=\frac{\int_{M}Ff\ dv_{\theta}\int_{M}FA_{\theta}^{-1}\tau(F)\ dv_{\theta}}{V^{2}_{F}}$.\\
\end{lemma}
{\it Proof:}
Recall that from our convention, we have that
$$A_{\theta_{F}}A_{\theta_{F}}^{-1}\tau_{F}f=\tau f-\frac{\int_{M}fdv_{\theta_{F}}}{V_{F}}.$$
But $A_{\theta_{F}}=\tau_{F}(F^{-1}A_{\theta})$. Hence, one has
$$A_{\theta}A_{\theta_{F}}^{-1}\tau_{F} f=\tau (Ff)-\frac{\tau(F)\int_{M}fdv_{\theta_{F}}}{V_{F}}.$$
Therefore,
\begin{align}
A^{-1}_{\theta_{F}}\tau_{F}f&=A_{\theta}^{-1}\tau (Ff)-\frac{\int_{M}fdv_{\theta_{F}}}{V_{F}}A_{\theta}^{-1}\tau(F)-\frac{1}{V_{F}}\int_{M}A_{\theta}^{-1}\tau(Ff)dv_{\theta_{F}}\notag\\
&\quad+\frac{\int_{M}fdv_{\theta_{F}}}{V_{F}^{2}}\int_{M}A_{\theta}^{-1}\tau F dv_{\theta_{F}}\notag.
\end{align}
\hfill$\Box$\\

We also recall here the scalar invariant related to the operator $P'_{\theta}$, namely, the $Q'$-curvature. Indeed, we set
 \begin{equation*}
  \label{qprime1}
  Q_\theta^\prime := 2\Delta_b R - 4 |A|^2 + R^2.
 \end{equation*}
Then for $w\in\mathcal{P}$ and $\hat\theta=e^w\theta$, we have
\begin{equation*}
 e^{2w} Q_{\hat{\theta}}^\prime = Q_\theta^\prime + P_\theta^\prime(w) + \frac{1}{2}P_\theta\left(w^2\right) .
\end{equation*}
In our case, we are more interested in the quantity $\bar{Q}'_{\theta}=\tau Q_{\theta}'$. For more information about the $Q'$-curvature we refer the reader to \cite{BG,CaYa}  and for problems related to prescribing the $\bar{Q}'_{\theta}$, we refer the reader to \cite{CaYa1,QQ,M}.
\begin{lemma}
Assume that $\ln(F)\in \P$, then we have
$$\tau_{F}Q'_{\theta_{F}}=\tau_{F}(F^{-1}Q'_{\theta})+\frac{1}{2}A_{\theta_{F}}\ln(F).$$
\end{lemma}
{\it Proof:}
Recall that under the conformal change $\theta\to \theta_{F}$, the $Q'$-curvature changes as follows:
$$\frac{1}{2}P'_{\theta}\ln(F)+Q'_{\theta}=Q'_{\theta_{F}}F+\frac{1}{8}P_{\theta}((\ln(F))^{2}).$$
Thus,
$$ Q'_{\theta_{F}}=\frac{F}{2}P'_{\theta}\ln(F)+FQ'_{\theta}-\frac{F}{8}P_{\theta}((\ln(F))^{2}).$$
Hence,
$$\tau_{F}Q'_{\theta_{F}}=\tau_{F}(F^{-1}Q'_{\theta})+\frac{1}{2}A_{\theta_{F}}\ln(F).$$
\hfill$\Box$

Define now the geometric mass as in \cite{S1,S2}, by
$$\mathcal{N}_{\theta}(x):=m_{\theta}(x)-\frac{\gamma_{3}}{2} A^{-1}_{\theta}\tau Q'_{\theta} (x).$$

A direct substitution then shows that if $\ln(F)\in \P$ then

\begin{align}
\mathcal{N}_{\theta_{F}}-\mathcal{N}_{\theta}&=\frac{\frac{\gamma_{3}}{2}\int_{M}Q'_{\theta}\ dv_{\theta}-2}{V_{F}}A_{\theta}^{-1}\tau F+\frac{1-\frac{\gamma_{3}}{2}\int_{M}Q'_{\theta}\ dv_{\theta}}{V_{F}^{2}}\int_{M}F A_{\theta}^{-1}\tau F \ dv_{\theta}\notag\\
&+\frac{\gamma_{3}}{4V_{F}}\int_{M}F\ln(F)\ dv_{\theta}-\frac{\gamma_{3}}{2V_{F}}\int_{M}FA_{\theta}^{-1}\tau Q'_{\theta}\ dv_{\theta}.\notag
\end{align}

In particulat, on the sphere $S^{3}$, since we have $\int_{S^{3}}Q'_{\theta}\ dv_{\theta}=16\pi^{2}$, we have 
\begin{proposition}
$$\mathcal{N}_{\theta_{F}}(S^{3})(x)-\mathcal{N}_{\theta_{0}}(S^{3})(x)=\frac{\gamma_{3}}{4V_{F}}\int_{M}F\ln(F)\ dv_{\theta}-\frac{1}{V_{F}^{2}}\int_{M}F A_{\theta}^{-1}\tau F \ dv_{\theta}\geq 0.$$
\end{proposition}

\end{document}